\documentclass[runningheads]{llncs}
\usepackage[T1]{fontenc}
\usepackage{graphicx}
\usepackage{algorithm}
\usepackage{algorithmicx}
\usepackage{algpseudocode}
\usepackage{amsmath}
\usepackage{amssymb}
\usepackage{optidef}
\usepackage{pgfplots}
\pgfplotsset{compat=1.18} 
\usepackage{tikz}
\usepackage{tikzscale}  
\usepackage{pgf-pie}   
\usepackage{tikz-3dplot}
\usetikzlibrary{patterns}
\usepackage{hyperref}
\usepgfplotslibrary{groupplots}
\usepackage{adjustbox}

\usepackage{subcaption}

\newcommand{\R}[1]{\mathbf{R}^{#1}}
\newcommand{\jac}[2]{\frac{\partial { #1  }}{ \partial { #2 }}}
\begin{document}
\title{On the Differentiability of the \\Primal-Dual Interior-Point Method}
%
%
\author{Kevin Tracy\and
Zachary Manchester}
\authorrunning{Tracy and Manchester}
%
\institute{The Robotics Institute, Carnegie Mellon University \\
\email{\texttt{\{ktracy, zacm\}}@cmu.edu}\\
}
\maketitle              
\begin{abstract}
Primal-Dual Interior-Point methods are capable of solving constrained convex optimization problems to tight tolerances in a fast and robust manner. 
The derivatives of the primal-dual solution with respect to the problem matrices can be computed using the implicit function theorem, enabling efficient differentiation of these optimizers for a fraction of the cost of the total solution time. 
In the presence of active inequality constraints, this technique is only capable of providing discontinuous subgradients that present a challenge to algorithms that rely on the smoothness of these derivatives. 
This paper presents a technique for relaxing primal-dual solutions with a logarithmic barrier to provide smooth derivatives near active inequality constraints, with the ability to specify a uniform and consistent amount of smoothing. 
We pair this with an efficient primal-dual interior-point algorithm for solving an always-feasible $\ell_1$-penalized variant of a convex quadratic program, eliminating the issues surrounding learning potentially infeasible problems.
This parallelizable and smoothly differentiable solver is demonstrated on a range of robotics tasks where smoothing is important. An open source implementation in JAX is available at \url{www.github.com/kevin-tracy/qpax}.
\end{abstract}
\section{Introduction}
Convex optimization has seen widespread use in modern robotics, where the guarantees of global optimality and polynomial time complexity have enabled algorithms that span control \cite{kuindersma2014a,blackmore2016}, state estimation \cite{varin2020,xinjilefu2014}, actuator allocation \cite{tracy2023d,kirchengast2018}, collision detection \cite{gilbert1994,tracy2023b}, and simulation \cite{anitescu2006,pang2021}. 

For years, practitioners have exploited domain-specific knowledge to craft convex optimization problems that enjoy fast and reliable convergence for both offline and online use. 
In the era of data-driven robotics, differentiable optimization has enabled automatic tuning or ``learning'' of optimization problems directly from data. Although this may not replace domain-specific knowledge, the ability to build complex tuneable functions with embedded optimization problems is well suited for a variety of tasks. 

The sensitivity analysis of linear programs has been studied for decades \cite{boyd2004}, but recent advances in differentiable convex optimization and parametrized convex layers in deep networks are gaining traction \cite{amos2019,agrawal2019a,agrawal2019}. This has resulted in convex modeling tools such as CVXPY layers that enable easy incorporation of differentiable convex optimization into common workflows \cite{diamond}. 

There are two main issues preventing the widespread adoption of differentiable optimization in robotics as it exists today: the subgradient problem and the infeasibility problem. For active inequality constraints in optimization problems, the derivatives propagated through these solvers are restricted to subgradients when there is no uniquely defined gradient. This results in routines getting ``stuck'' near these inequalities due to their nonsmooth nature. The second issue arises when infeasible problem instances are created during auto-tuning. While infeasibility detection is commonplace in convex optimization, detecting an infeasible problem does not address the root cause of the infeasibility, nor provide informative gradient information to encourage feasibility. 

In this paper, we propose solutions to both of these problems. Our contributions include:
\begin{itemize}
    \item A rigorous method for returning unique and smoothed derivatives of convex optimization problems through the use of a relaxed logarithmic barrier 
    \item A primal-dual interior-point algorithm for solving an always-feasible ``elastic'' quadratic program with minimal computational overhead
\end{itemize}
Together, these two advances are presented in an open source software package written in JAX and demonstrated on the relevant robotics tasks of contact dynamics and collision detection. 
%
%
%
%
\section{Background}
This section introduces the notation for a standard form of the convex Quadratic Program (QP), common solution methods, and extensions to differentiable optimization.
\subsection{Quadratic Programming}
In this paper, we focus on the convex QP as a fundamental problem specification. Modern algorithms such as the Primal-Dual Interior-Point (PDIP) method solve these problems globally in a fast and efficient manner \cite{mattingley2012}. A standard form quadratic program and its equivalent with a slack variable are as follows:

\begin{tabular}{c c c}
    \begin{minipage}[t]{0.4\textwidth}
        \centering 
        \begin{mini}
            {x}{ \frac{1}{2}x^TQx + q^Tx }{\label{qp_standard_form}}{}
            \addConstraint{Ax}{=b}
            \addConstraint{Gx}{\leq h,}
        \end{mini}
    \end{minipage}
    &
    \begin{minipage}[t]{0.1\textwidth}
        \centering
        \vspace{40pt}
        $\Rightarrow$ 
    \end{minipage}
    \begin{minipage}[t]{0.4\textwidth}
        \centering
        \begin{mini}
            {x, s}{ \frac{1}{2}x^TQx + q^Tx }{\label{qp_standard_form_2}}{}
            \addConstraint{Ax}{=b}
            \addConstraint{Gx + s}{= h}
            \addConstraint{s}{\geq 0,}
        \end{mini}
    \end{minipage}
\end{tabular}

with a primal variable $x \in \R{n}$, cost terms $Q \in \mathbf{S}_+^n$ and $q \in \R{n}$, equality constraints described with $A \in \R{m \times n}$ and $b \in \R{m}$, and inequality constraints with $G \in \R{p \times n}$ and $h \in \R{p}$. Dual variables are introduced to enforce the constraints, with $y \in \R{m}$ associated with the equality constraint, and $z \in \R{p}$ with the inequality constraints \cite{boyd2004}. The slack variable $s \in \R{p}$ is introduced for algorithmic simplicity.
The Lagrangian for this problem is then
\begin{align}
    \mathcal{L}(x,s,z,y) &= \frac{1}{2}x^TQx + q^Tx + y^T(Ax - b) + z^T(Gx -h),
\end{align}
resulting in the following KKT conditions for optimality:
\begin{align}
    Qx + q + G^Tz + A^Ty &= 0, \label{qp:kkt:stat}\\ 
    z \odot  s &= 0, \label{qp:kkt:compl} \\ 
    Gx + s &= h, \label{qp:kkt:pfeas1}\\ 
    Ax &= b,  \label{qp:kkt:pfeas2} \\ 
    s &\geq 0, \label{qp:kkt:pfeas3}\\ 
    z &\geq 0, \label{qp:kkt:dfeas}
\end{align}
where $\odot$ denotes elementwise multiplication. A primal-dual solution $(x^*, s^*, y^*, z^*)$ is globally optimal if it satisfies \eqref{qp:kkt:stat}-\eqref{qp:kkt:dfeas}.
\subsection{Primal-Dual Interior-Point Methods}
PDIP methods solve \eqref{qp_standard_form_2} by treating a modified version of the system of equations in \eqref{qp:kkt:stat}-\eqref{qp:kkt:pfeas2} as a root finding problem, and then using Newton's method to find a solution while restricting $(s, z) > 0$. As a result, the majority of the computation time is spent solving linear systems of the following form:
\begin{align}
    \begin{bmatrix}
Q & 0 & G^{T} & A^{T} \\
0 & D(z) & D(s) & 0 \\
G & I & 0 & 0 \\
A & 0 & 0 & 0
\end{bmatrix} \begin{bmatrix}
\Delta x\\
\Delta s \\
\Delta z \\
\Delta y
\end{bmatrix}=\begin{bmatrix}
u_1 \\ u_2 \\ u_3 \\ u_4
\end{bmatrix}, \label{normal_pdip_ls}
\end{align}
where $D(\cdot)$ denotes the diagonal matrix constructor from a vector. Using block reduction techniques, the linear system in \eqref{normal_pdip_ls} can be efficiently solved with Alg. \eqref{alg:normal_pdip_ls_solve}. This technique for solving linear systems of this form is useful for both the solving and differentiation of the PDIP method. 

After the step directions are computed with the linear system, a linesearch is used to ensure the nonnegativity of $(s, z)$ .  For an arbitrary variable $v$ and step direction $\Delta v$, a linesearch that solves for the largest $\alpha \leq 1$ that keeps $v + \Delta v\geq 0$ is solved in closed form:
%
\begin{align}
    \operatorname{linesearch}(v, \Delta v) &= \min \bigg( 1,\, \min_{i:\Delta v_i < 0} -\frac{v_i}{\Delta v_i} \bigg) .\label{sec:background:linesearch}
\end{align}
This linesearch is performed for both $s$ and $z$, and the step length is simply the minimum of the two. For a complete algorithmic specification and implementation details, we refer the reader to \cite{mattingley2012} for a PDIP method for which Alg. \eqref{alg:normal_pdip_ls_solve} can be used to solve the linear systems. 
\begin{algorithm}
\caption{PDIP Linear System Solver}\label{alg:normal_pdip_ls_solve}
\begin{algorithmic}[1]
\State \textbf{function} $\operatorname{solve\_kkt}(u_1, u_2, u_3, u_4)$ \Comment{assume access to $Q,G,A,s,z$}
\State $P \gets D(s \oslash  z)$ 
\State $H \gets Q + G^TP^{-1}G$ \Comment{cacheable Cholesky decomposition}
\State $F \gets A H^{-1} A^T$  \Comment{cacheable Cholesky decomposition}
\State $r_2 \gets u_3 - u_2 \oslash z$ 
\State $p_1 \gets u_1 + G^T P^{-1}r_2$ 
\State $\Delta y \gets F^{-1}(A H^{-1}p_1 - u_4)$
\State $\Delta x \gets H^{-1}(p_1 - A^T\Delta y)$ 
\State $\Delta s \gets u_3 - G \Delta x$ 
\State $\Delta z \gets (u_2 - z \odot \Delta s) \oslash s$
\State \textbf{return} $(\Delta x, \Delta s, \Delta z, \Delta y)$
\end{algorithmic}
\end{algorithm}
\subsection{Differentiable Optimization}
Modern automatic differentiation tools have made forming derivatives through complex functions easier than ever. When one of these functions includes an iterative routine, unrolling the iterations into a sequential computational graph and proceeding with differentiation is not always possible. The first concern is nonsmooth operations or logical branching in the iterations that can stop the ``flow'' of the derivatives through the routine, and the second concern is the rapidly decaying numerical precision inherent in differentiating an iterative process.

To address both of these shortcomings, iterative routines such as numerical optimizers are differentiated with methods that do not require propagating derivatives through the iterations themselves. Instead, the optimization problem is solved as normal, and the derivatives are then constructed directly from the solution to the problem.
\subsubsection{Implicit Function Theorem}
When an iterative routine can be interpreted as finding a root or an equilibrium point to an implicit function, the implicit function theorem can be used to form the derivatives of interest. Given variables $w \in \R{a}$ and parameters $\theta \in \R{b}$, an implicit function is defined as
\begin{align}
    r(w^*, \theta) = 0,
\end{align}
at an equilibrium point $w^*$. This implicit function can be linearized about this point resulting in the following first-order Taylor series:
\begin{align}
    \jac{r}{w} \delta w + \jac{r}{\theta} \delta \theta = 0 ,
\end{align}
which can simply be re-arranged to solve for
\begin{align}
    \jac{w}{\theta} = -\bigg(\jac{r}{w} \bigg)^{-1} \jac{r}{\theta}. \label{ift}
\end{align}
The implicit function theorem enables the differentiation of routines that solve for equilibrium points without the need to unroll and differentiate the iterations themselves. 

By treating the KKT conditions from \eqref{qp:kkt:stat}-\eqref{qp:kkt:pfeas2} as a residual function of the primal-dual solution $(x, s, z, y)$ and problem parameters $\theta$, the implicit function theorem is used to form the derivatives of the optimizer without unrolling the iterations. The linear system from \eqref{ift} applied to this residual function is of the same form as that in Alg. \eqref{alg:normal_pdip_ls_solve}, enabling fast and easy computation of derivatives using a framework already available in the solver.
%
%
\subsubsection{Efficient Computation of Gradients}
Using the implicit function theorem to compute the Jacobians of the primal-dual solution with respect to the problem parameters results in the need to directly form potentially large Jacobians. In many cases, it is not these specific Jacobians that are of interest, but rather the left matrix-vector product with a backward pass vector.

For a loss function of the primal variable $\ell(x)$ that takes as input the optimal primal variable from \eqref{qp_standard_form}, reverse-mode automatic differentiation will have constructed $\nabla_x \ell$ by the time it comes to the QP solver in the backward pass of the computational graph. Instead of forming the Jacobians of the primal variable with respect to each of the problem matrices directly, it is instead desirable to form the 
left matrix-vector products $\jac{\ell}{x} \jac{x}{\square}$,
where $\square$ simply denotes any of the problem matrices \cite{amos2017}.

As shown in Alg. \eqref{alg:opt_net}, the gradients of this loss function $\nabla_\square \ell$ with respect to the problem parameters of the QP can be computed by once again utilizing our interior-point linear system solver.
\begin{algorithm} 
\begin{algorithmic}[1]
    \caption{Computing Gradients Through a QP}\label{alg:opt_net}
        \State \textbf{function} $\operatorname{compute\_qp\_grads}(Q, q, A, b, G, h, x, s, z, y, \nabla_x \ell)$ 
        \State $dx, ds, d\tilde{z}, dy \gets \operatorname{solve\_kkt}(-\nabla_x \ell, 0, 0, 0)$ \Comment{compute differentials with kkt system}
        \State $dz = d\tilde{z} \oslash z$ 
        \State $\nabla_Q \ell \gets \big((dx) x^T + x (dx)^T)\big)/2$
        \State $\nabla_q \ell \gets dx$
        \State $\nabla_A \ell \gets (dy)x^T + y (dx)^T$ 
        \State $\nabla_b \ell \gets -dy$
        \State $\nabla_G \ell \gets z \odot \big((dz)x^T + z(dx)^T\big)$
        \State $\nabla_h \ell \gets -z \odot dz$
    \State \textbf{return} $\nabla_Q \ell, \nabla_q \ell, \nabla_A \ell, \nabla_b \ell, \nabla_G \ell, \nabla_h \ell$ 
\end{algorithmic}
\end{algorithm}
\section{Logarithmic Barrier Smoothing}
\begin{figure}[t!]
    \centering
    \hspace{0cm}
    \includegraphics[width=.4\linewidth]{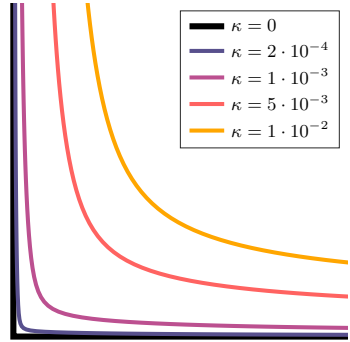}
    \caption{A sharp corner in a square as smoothed with the logarithmic barrier at varying central path parameters $\kappa$. As $\kappa \rightarrow 0$, the corner becomes more pronounced until it assumes a true $90^\circ$ corner at $\kappa = 0$. The logarithmic barrier effectively smooths out any sharp corners of the feasible set enabling smooth differentiation in the presence of such discontinuities.}
    \label{fig:log_barrier}
\end{figure}
Primal interior-point methods work by replacing constraints of the form $h(x) \leq 0$ with a logarithmic barrier penalty in the form of $\phi(x) = - \kappa \sum \log(-h(x))$, where $\kappa \in \R{+}$ is referred to as the central path parameter. This barrier function is a smooth approximation of the indicator function, where feasible values of $x$ result in no penalty, and infeasibility results in an infinite penalty \cite{boyd2004}. By solving a sequence of unconstrained problems as $\kappa \rightarrow 0$, the logarithmic barrier becomes a closer and closer approximation of the indicator function until acceptable convergence is achieved. This barrier function is only defined for feasible values of $x$, hence the name ``interior-point''. 

For a quadratic program in standard form \eqref{qp_standard_form_2}, the optimality conditions for a barrier subproblem given a central path parameter $\kappa$ are 
almost identical to the original KKT conditions \eqref{qp:kkt:stat}-\eqref{qp:kkt:dfeas} with the exception of the complementarity condition \eqref{qp:kkt:compl} replaced with $z \odot s - \kappa = 0$.  This relaxed complementarity condition allows for a certain amount of smoothing over the feasible set, where larger values of $\kappa$ have a stronger smoothing effect. The optimality conditions for this barrier subproblem are referred to as the perturbed or relaxed KKT conditions. 
\subsection{Relaxing Primal-Dual Solutions}
\begin{algorithm}
\caption{Relaxing a Quadratic Program}\label{alg:relax_qp}
\begin{algorithmic}[1]
\State \textbf{function} $\operatorname{relax\_qp}(Q,q,A, b, G,h, x, s, z, y, \kappa)$ 
\For{$i \gets 1:\texttt{max\_iters}$} 
\State \texttt{/* evaluate relaxed KKT conditions and check convergence*/}
\State $r_1 \gets Qx + q  + A^Ty + G^Tz$ 
\State $r_2 \gets s \odot z - \kappa$ \Comment{relaxed complementarity}
\State $r_3 \gets Gx + s - h$
\State $r_4 \gets Ax - b$
\If{$\|(r_1, r_2, r_3, r_4, r_5, r_6)\|_\infty < \texttt{tol}$} 
\State \textbf{return:} $x, t, s_1, s_2, z_1, z_2$
\EndIf 
\State 
\State \texttt{/* calculate and take Newton step */}
\State $\Delta x, \Delta s, \Delta z, \Delta y \gets \operatorname{solve\_kkt}(-r_1, -r_2, -r_3, -r_4)$ \Comment{Alg. \eqref{alg:normal_pdip_ls_solve}}
\State $\alpha \gets 0.98 \min\big(\operatorname{linesearch}(s, \Delta s),\, \operatorname{linesearch}(z, \Delta z)\big)$ \Comment{\eqref{sec:background:linesearch}}
\State $(x, s, z, y) \gets (x, s, z, y) + \alpha(\Delta x, \Delta s, \Delta z, \Delta y)$
\EndFor
\end{algorithmic}
\end{algorithm}
When differentiating a quadratic program near an active inequality constraint, the implicit function theorem can produce subgradients that are potentially uninformative \cite{agrawal2019a}. In order to ensure smooth and continuous gradients in the presence of sharp corners in the feasible set, proposed here is a relaxation method that exploits the smoothing of the logarithmic barrier to effectively round out any nonsmoothness. Specifically, the idea is to take a primal-dual solution that is optimal to some low $\kappa_\text{low}$, and relax it to a specified $\kappa_\text{high}$.

In PDIP methods, much care is taken to solve the perturbed KKT conditions for decreasing $\kappa \rightarrow 0$. The most common strategy in PDIP methods is a Mehrotra predictor-corrector method that adaptively updates the target $\kappa$ until it is below the convergence criteria \cite{mehrotra1992}. While solving this system for a sequence where $\kappa \rightarrow 0$ is challenging, going the other direction from $\kappa_\text{low} \rightarrow \kappa_\text{high}$ is actually quite trivial. For this case, standard Newton steps on the perturbed KKT conditions with a linesearch to ensure $(s,z) > 0$ is able to converge to $\kappa_\text{high}$ in only a few steps. This algorithm is shown in Alg. \eqref{alg:relax_qp}, where once again the linear system solver from Alg. \eqref{alg:normal_pdip_ls_solve} is used. 

Once the primal-dual solution to the quadratic program has been relaxed, Alg. \eqref{alg:opt_net} is used to calculate gradients that benefit from the logarithmic barrier smoothing. The full sequence for the solving, relaxation, and differentiation of a quadratic problem is as follows:
\begin{enumerate}
    \item Solve the quadratic program to a specified tolerance and return the solution
    \item Relax the primal-dual solution to a target $\kappa$ 
    \item Form the derivatives of interest at the relaxed primal-dual solution
\end{enumerate}
This means that the solver itself can return high quality solutions to tight tolerances while still returning smooth gradients evaluated at the relaxed solution. In modern automatic differentiation frameworks, this sequence is written into a custom forward and backward pass through the function with limited overhead. 

The choice of a target $\kappa$ is left to the specifics of the problem at hand. In some scenarios, only a little bit of smoothing is required, making a lower $\kappa$ appropriate. For very sharp corners in the feasible set (like the tip of a triangle), a larger value of $\kappa$ can be used to provide even more smoothing. In either case, computing the derivatives of the solver with this technique allows both the tolerance of the solver and the relaxed $\kappa$ to be specified exactly and independently. 
\section{Elastic Quadratic Program Solver}
While the convexity of a quadratic program guarantees a globally optimal solution when one is available, there is generally no guarantee of feasibilty. This infeasibility occurs when the quadratic program has a set of constraints that are impossible to satisfy, something that can easily happen in practice if the constraint matrices are learned. In the event of an infeasible problem, standard PDIP methods are unable to return a useful solution. 

Traditionally, infeasibility in convex optimization problems can be handled with a homogenous self-dual embedding that allows for the computation of a certificate of infeasibility \cite{domahidi2013,vandenberghe,stellato}. This ensures that, even in the event of an infeasible problem, the solver can ``gracefully'' fail in a way that simply returns certification that there are no solutions that satisfy the constraints. While this approach is useful for determining if a given problem is feasible, ideally we would set up problems such that infeasibility is not possible. 

Given a convex quadratic program in a standard inequality-only form, 
\begin{mini}
    {x}{ \frac{1}{2}x^TQx + q^Tx }{\label{qp_ineq_only}}{}
    \addConstraint{Gx}{\leq h,}
\end{mini}
it is possible there is no $x$ that satisfies $Gx\leq h$. In order to convert the optimization problem in \eqref{qp_ineq_only} into one in which there is always a solution, the hard constraints are converted into penalties. An $\ell_1$-penalty on the constraint violation is chosen because the $\ell_1$-norm encourages sparsity in the constraint violation, translating into a penalty that encourages the solver to satisfy as many of the constraints as possible. This is a common technique in nonlinear programming \cite{nocedal2006} for the handling of infeasible subproblems. This ``elastic'' mode, as it is referred to in the solver SNOPT \cite{gill2005}, is a highly effective method to guarantee that a problem always has a solution without sacrificing the quality and utility of the solution. The problem from \eqref{qp_ineq_only} is converted into its elastic form as follows:
\begin{mini}
    {x}{ \frac{1}{2}x^TQx + q^Tx + \|\rho \odot \max(0, Gx-h)\|_1,}{\label{qp_elastic_unconstrained}}{}
\end{mini}
such that feasible values of $x$ do not contribute to the cost function but infeasibility is penalized with $\rho \in \R{p}$. While \eqref{qp_elastic_unconstrained} is an unconstrained convex optimization problem, it is nonsmooth and does not pair well with PDIP methods. We therefore reformulate it as:
\begin{mini}
    {x, t}{ \frac{1}{2}x^TQx + q^Tx + \rho^T t}{\label{qp_elastic_constrained}}{}
    \addConstraint{Gx - h}{\leq t}
    \addConstraint{t}{\geq 0,}
\end{mini}
where $t \in \R{p}$ is a slack variable containing the constraint violation, and a simple linear cost term is used to recover the $\ell_1$-penalty from \eqref{qp_elastic_unconstrained}.

Solving the elastic quadratic program in \eqref{qp_elastic_constrained} with a PDIP method normally comes at an added computational expense since we have increased the number of primal-dual variables from $n + 2p$ to $n + 5p$. To avoid the cubic complexity in the increase primal-dual dimension, we introduce a custom algorithm for solving these problems that exploits the sparsity of the constraints such that the time to solve the elastic version of the problem is only a slight (5--20\%) increase compared to the time to solve the original problem.

As before, the PDIP method for solving \eqref{qp_elastic_constrained} is dominated by the factorization and solving of linear systems in the following form:
    \begin{align}
        \begin{bmatrix} 
            Q & 0 & 0   & 0   & 0   & G^T \\ 
            0 & 0 & 0   & 0   & -I   & -I \\ 
            0 & 0 & Z_1 & 0   & S_1 & 0 \\ 
            0 & 0 & 0   & Z_2 & 0   & S_2 \\ 
            0 & -I & I & 0 & 0 & 0 \\ 
            G & -I & 0 & I & 0 & 0 
        \end{bmatrix} \begin{bmatrix}
            \Delta x \\ \Delta t \\ \Delta s_1 \\ \Delta s_2 \\ \Delta z_1 \\ \Delta z_2
        \end{bmatrix} &= \begin{bmatrix} r_1 \\ r_2 \\ r_3 \\ r_4 \\ r_5 \\ r_6 \end{bmatrix}.
    \end{align}
    As shown in Alg. \eqref{alg:elastic_kkt}, this linear system can be solved with block-wise elimination where the only matrix factorization required is that of a positive definite matrix the size of the primal variable $x$. This routine is used in the full PDIP algorithm for the elastic QP as shown in Alg. \eqref{alg:elastic_pdip}. 

    The elastic QP is fully differentiable in the same way the original QP is. To do this, Alg. \eqref{alg:relax_elastic_qp} is used to take a primal-dual solution and relax it to a specified $\kappa$, and  Alg. \eqref{alg:elastic_kkt} for solving the linear systems is re-used for fast and efficient relaxation of the elastic problem. From this, the gradients of a downstream loss function with respect to the problem parameters can be constructed with Alg. \eqref{alg:opt_net_elastic}, again using the same linear system solver. 

    The solving, relaxation, and differentiation of the elastic mode QP are all only a slight increase in computational complexity compared to the original QP, with the benefit of guaranteed feasibility. This enables the inclusion of always-feasible quadratic programming in learned pipelines where feasibility cannot be guaranteed by construction.
%
%
\begin{algorithm}
\caption{PDIP Method for Elastic Quadratic Programs}\label{alg:elastic_pdip}
\begin{algorithmic}[1]
\State \textbf{function} $\operatorname{solve\_qp\_elastic}(Q,q, G,h, \rho)$ 
\State $x, s_1, s_2, z_1, z_2 \gets \operatorname{initialize}(Q,q, G,h, \rho)$ \Comment{\ref{sec:elastic_init}}
\For{$i \gets 1:\texttt{max\_iters}$} 
\State \texttt{/* evaluate KKT conditions and check convergence*/}
\State $r_1 \gets Qx + q  + G^Tz_2$ 
\State $r_2 \gets -z_1 - z_2 + \rho$ 
\State $r_3 \gets s_1 \odot z_1$ 
\State $r_4 \gets s_2 \odot z_2$ 
\State $r_5 \gets -t + s_1 $ 
\State $r_6 \gets Gx - t + s_2 - h$
\If{$\|(r_1, r_2, r_3, r_4, r_5, r_6)\|_\infty < \texttt{tol}$} 
\State \textbf{return:} $x, t, s_1, s_2, z_1, z_2$
\EndIf 
\State 
\State \texttt{/* calculate affine step direction Alg. \eqref{alg:elastic_kkt}*/}
\State $\Delta x^\text{a}, \Delta t^\text{a}, \Delta s_1^\text{a}, \Delta s_2^\text{a}, \Delta z_1^\text{a}, \Delta z_2^\text{a}\gets \operatorname{elastic\_kkt}(-r_1, -r_2, -r_3, -r_4, -r_5, -r_6)$ 
\State $\Delta s^\text{a},\, \Delta z^\text{a} \gets (\Delta s_1^\text{a}, \Delta s_2^\text{a}), \,(\Delta z_1^\text{a}, \Delta z_2^\text{a})$
\State $s,\, z \gets (s_1,s_2), \, (z_1, z_2)$
\State \texttt{/* calculate centering-plus-corrector step direction */}
\State $\alpha^\text{a} = \min( \operatorname{linesearch}(s, \Delta s^\text{a}),\, \operatorname{linesearch}(z, \Delta z^\text{a}) )$ \Comment{\eqref{sec:background:linesearch}}
\State $\mu \gets s^Tz/\operatorname{len}(s)$
\State $\sigma \gets [(s + \alpha^\text{a} \Delta s^\text{a})^T(z + \alpha^\text{a}\Delta z^\text{a})/(s^Tz)]^3$
\State $r_3 \gets r_3 - \sigma \mu \mathbf{1} + \Delta s_1^\text{a} \odot \Delta z_1^\text{a}$
\State $r_4 \gets r_4 - \sigma \mu \mathbf{1} + \Delta s_2^\text{a} \odot \Delta z_2^\text{a}$
\State $\Delta x, \Delta t, \Delta s_1, \Delta s_2, \Delta z_1, \Delta z_2 \gets \operatorname{elastic\_kkt}(-r_1, -r_2, -r_3, -r_4, -r_5, -r_6)$ 
\State 
\State \texttt{/* update with linesearch */}
\State $\alpha \gets 0.98 \min(\operatorname{linesearch}(s, \Delta s),\, \operatorname{linesearch}(z, \Delta z))$ 
\State $(x, t, s_1, s_2, z_1, z_2) \gets (x, t, s_1, s_2, z_1, z_2)  + \alpha(\Delta x, \Delta t, \Delta s_1, \Delta s_2,  \Delta z_1, \Delta z_2)$
\EndFor
\end{algorithmic}
\end{algorithm}
\begin{algorithm} 
\begin{algorithmic}[1]
    \caption{Solve Elastic KKT Linear System}\label{alg:elastic_kkt}
        \State \textbf{function} $\operatorname{elastic\_kkt}(r_1, r_2, r_3, r_4, r_5, r_6)$ \Comment{assume access to $Q,G,A,s_1,s_2,z_1,z_2$}
        \State $w_1 \gets r_3 \oslash z_1$
        \State $w_2 \gets r_4 \oslash z_2$
        \State $p_1 \gets r_5 - r_6 + w_2 - w_1 - (s_1 \odot r_2) \oslash z_1$
        \State $A_3 \gets \operatorname{diag}(a_1 + a_2)$
        \State $\Delta x \gets (Q + G^TA_3^{-1}G)^{-1}(r_1 - G^TA_3^{-1}p_1) $ \Comment{Cholesky factorization (cacheable)}
        \State  $\Delta z_2 \gets A_3^{-1}(p_1 + G \Delta x) $ 
        \State  $\Delta z_1 \gets -r2 - \Delta z_2 $
        \State  $\Delta s_1 \gets (r_3 - s_1 \odot \Delta z_1) / z_1  $ 
        \State  $\Delta s_2 \gets (r_4 - s_2 \odot \Delta z_2) / z_2 $
        \State  $\Delta t \gets \Delta s_1 - r_5 $
    \State \textbf{return} $\Delta x,\, \Delta s_1,\, \Delta s_2,\, \Delta z_1,\, \Delta z_2$ 
\end{algorithmic}
\end{algorithm}
\subsection{Elastic Initialization}\label{sec:elastic_init}
In order to initialize the PDIP method shown in Alg. \eqref{alg:elastic_pdip}, the only requirement is that $s,z > 0$. In practice, a more advanced initialization technique can both reduce the number of iterations required for convergence and dramatically improve the robustness of the solver. The initialization from \cite{vandenberghe} and \cite{mattingley2012} is adapted for use in the elastic case. 

First, the solution to the following problem is computed analytically:
\begin{mini}
    {x, t, s_1, s_2}{ \frac{1}{2}x^TQx + q^Tx + \frac{1}{2}(s_1^Ts_1 + s_2^Ts_2)  + \rho^T t}{\label{pdip_init_qp_problem}}{}
    \addConstraint{s_1 - t}{=0}
    \addConstraint{Gx - t + s_2}{= h,}
\end{mini}
where the primal and dual solutions are the solution to the linear system,
\begin{align}
\begin{bmatrix} 
    Q & 0    & 0   & G^T \\ 
    0 & 0    & -I   & -I \\ 
    0 & -I & -I & 0    \\ 
    G & -I & 0 & -I  
\end{bmatrix} \begin{bmatrix}
    \Delta x \\ \Delta t \\ \Delta z_1 \\ \Delta z_2
\end{bmatrix} &= \begin{bmatrix} -q \\ \rho \\ 0 \\ h \end{bmatrix},
\end{align}
which can be solved with a dense block reduction:
\begin{align}
    x &= (Q + \frac{1}{2}G^TG)^{-1}\big(-q - \frac{1}{2}G^T(\rho - h)\big) ,\\ 
    z_2 &= \frac{1}{2}(Gx + \rho - h), \\ 
    z_1 &= \rho - z_2, \\ 
    t &= -z_1 .
\end{align}
This method only requires one positive definite matrix the size of $x$. From here, we stack $z = (z_1, z_2)$ and initialize $s = (s_1, s_2)$ with
\begin{align}
    \alpha_p &= - \min(-z), \\ 
    s &= \begin{cases} -z, & \alpha_p < 0 \\ 
                       -z + (1 + \alpha_p \mathbf{1}), & \alpha_p \geq 0
                       \end{cases}. 
\end{align}
The dual variable for the inequality constraint is then initialized with 
\begin{align}
    \alpha_d &= - \min(-z), \\ 
    z &= \begin{cases} z + (1 + \alpha_d \mathbf{1}), & \alpha_d \geq 0 \\ 
                       z, & \alpha_d < 0
                       \end{cases}, 
\end{align}
finishing the initialization of the primal and dual variables.

\begin{algorithm}
\caption{Relaxing an Elastic Quadratic Program}\label{alg:relax_elastic_qp}
\begin{algorithmic}[1]
\State \textbf{function} $\operatorname{relax\_qp\_elastic}(Q,q,A, b, G,h, \rho, x, t, s_1, s_2, z_1, z_2, \kappa)$ 
\For{$i \gets 1:\texttt{max\_iters}$} 
\State \texttt{/* evaluate KKT conditions and check convergence*/}
\State $r_1 \gets Qx + q  + G^Tz_2$ 
\State $r_2 \gets -z_1 - z_2 + \rho$ 
\State $r_3 \gets s_1 \odot z_1 - \kappa$ \Comment{relaxed complementarity} 
\State $r_4 \gets s_2 \odot z_2 - \kappa$ \Comment{relaxed complementarity} 
\State $r_5 \gets -t + s_1 $ 
\State $r_6 \gets Gx - t + s_2 - h$
\If{$\|(r_1, r_2, r_3, r_4, r_5, r_6)\|_\infty < \texttt{tol}$} 
\State \textbf{return:} $x, t, s_1, s_2, z_1, z_2$
\EndIf 
\State 
\State \texttt{/* calculate and take Newton step*/}
\State $\Delta x, \Delta t, \Delta s_1, \Delta s_2, \Delta z_1, \Delta z_2\gets \operatorname{elastic\_kkt}(-r_1, -r_2, -r_3, -r_4, -r_5, -r_6)$ 
\State $\Delta s,\, \Delta z\gets (\Delta s_1, \Delta s_2), \,(\Delta z_1, \Delta z_2)$
\State $s,\, z \gets (s_1,s_2), \, (z_1, z_2)$
\State $\alpha^\text{a} = 0.98 \min( \operatorname{linesearch}(s, \Delta s),\, \operatorname{linesearch}(z, \Delta z) )$ \Comment{\eqref{sec:background:linesearch}}
\State $(x, t, s_1, s_2, z_1, z_2) \gets (x, t, s_1, s_2, z_1, z_2)  + \alpha(\Delta x, \Delta t, \Delta s_1, \Delta s_2,  \Delta z_1, \Delta z_2)$
\EndFor
\end{algorithmic}
\end{algorithm}
\begin{algorithm} 
\begin{algorithmic}[1]
    \caption{Computing Gradients Through an Elastic QP}\label{alg:opt_net_elastic}
        \State \textbf{function} $\operatorname{compute\_elastic\_qp\_grads}(Q, q, G, h, x, t,s_1,s_2,z_1,z_2,\nabla_x \ell)$ 
        \State \texttt{/* compute differentials with same linear system */}
        \State $dx, dt, ds_1, ds_2, d\tilde{z}_1, d\tilde{z}_2 \gets \operatorname{elastic\_kkt}(-\nabla_x \ell, 0, 0, 0, 0, 0)$
        \State $p, s, z \gets (x, t), (s_1, s_2), (z_1, z_2)$ 
        \State $dp, ds, dz \gets (dx, dt), (ds_1, ds_2), (d\tilde{z}_1, d\tilde{z}_2) \oslash (z_1, z_2)$
        \State \texttt{/* create indices (one-based) and form gradients from differentials */}
        \State $ix \gets 1:\operatorname{len}(q)$ 
        \State $is \gets (\operatorname{len}(q) + 1):(\operatorname{len}(q) + \operatorname{len}(h))$
        \State $\nabla_Q \ell \gets \big[ \big((dp) p^T + p (dp)^T)\big)/2 \big]_{ix,ix}$
        \State $\nabla_q \ell \gets [dp]_{ix}$
        \State $\nabla_G \ell \gets \big[ z \odot \big((dz)p^T + z(dp)^T\big) \big]_{is, ix}$
        \State $\nabla_h \ell \gets [-z \odot dz]_{is}$
    \State \textbf{return} $\nabla_Q \ell, \nabla_q \ell, \nabla_A \ell, \nabla_b \ell, \nabla_G \ell, \nabla_h \ell$ 
\end{algorithmic}
\end{algorithm}
\section{Numerical Experiments}
The utility of the proposed approaches are demonstrated on the common optimization-based robotics tasks of contact mechanics simulation and collision detection. In each of these tasks, a QP makes up a core part of the algorithm and smooth differentiation through the inherent nonsmoothness proves useful. 
\subsection{Contact Mechanics}
    \begin{figure}[t!]
        \centering
        \hspace{1cm}
        \includegraphics[width=.6\linewidth]{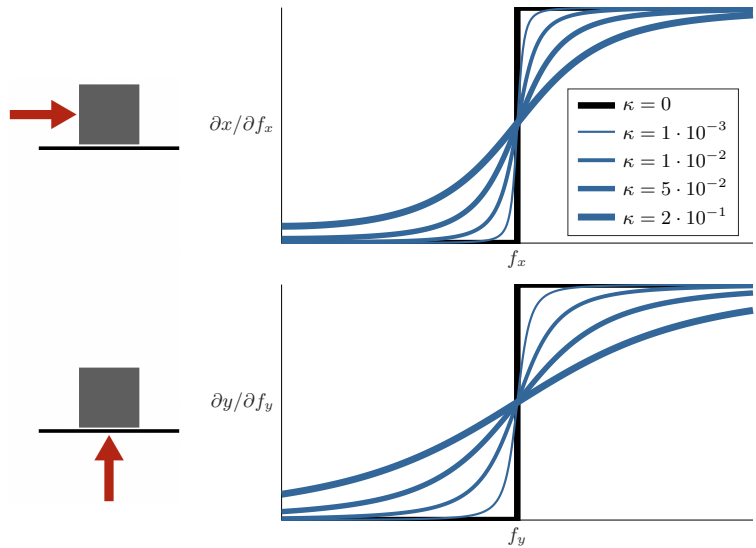}
        \caption{Contact dynamics for a two-dimensional block as modeled with a quadratic program. When a horizontal force $f_x$ is applied to the block, it must overcome the friction with the ground before it moves, and when a vertical force $f_y$ is applied, it must overcome gravity. Despite these discontinuities in the dynamics, the relaxed gradients from the differentiable quadratic program solver are able to provide smooth and continuous derivative information before and after the block begins to move.}
        \label{fig:dojo}
    \end{figure}
    \begin{figure}
        \centering
        \includegraphics[width=.9\linewidth]{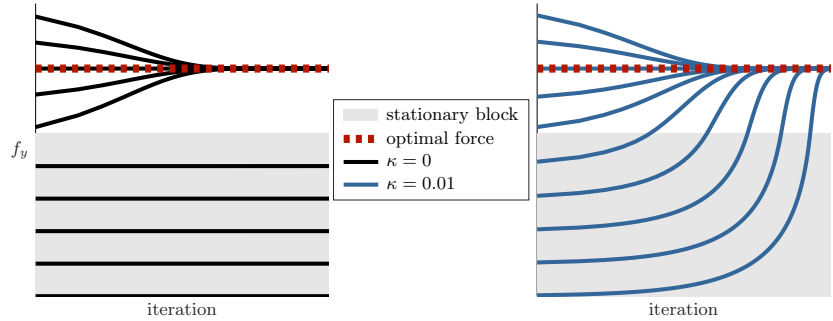}
        \caption{Using the gradients from the block pushing example in Fig. \ref{fig:dojo}, a vertical force is optimized to accelerate the block to a set value from multiple different initial force values. When the force is below the threshold for the block to move and $\kappa = 0$, the gradient is zero and the optimizer fails to make progress. In the case of $\kappa=0.01$, even before the block moves there is gradient information that pushes the optimizer to converge on the optimal force regardless of the initial force value.}
        \label{fig:dojo_convergence}
    \end{figure}
For a block at rest on a table with both gravity and friction, the nonsmooth contact dynamics can be represented as the solution to a convex quadratic program. For a full treatment of optimization-based dynamics, the reader is referred to \cite{anitescu2006} and \cite{howell2022}.

In the simplest case, a block is stationary until it is acted upon by a force that exceeds the static frictional force in the horizontal direction, or the gravitational force in the vertical direction. Until the applied forces exceed these two thresholds, the block does not move. This exercise is demonstrated in Fig. \ref{fig:dojo}, where forces are applied in each of the two directions and the true derivatives of these dynamics at $\kappa=0$ have a discontinuity as soon as the block begins to move. When these derivatives are taken with a relaxed $\kappa > 0$, the discontinuous derivative is smoothed, allowing for an informative gradient about the impending motion of the block before it moves. 

This derivative information is used in an optimization routine in Fig. \ref{fig:dojo_convergence}, where an optimizer attempts to solve for an applied force that produces the desired motion from the block. This optimizer is initialized at multiple different force values, and uses either exact gradients with $\kappa=0$, or smoothed gradients with $\kappa=0.01$. In the case of exact gradients, for initial forces where the block does not move, the gradient is zero and the solver is unable to find a descent direction. When $\kappa$ is relaxed, the smooth gradients provide information about the motion of the block even before the block itself begins to move, allowing the optimizer to converge on the true solution for each initial force. This is a simple yet expressive demonstration of the impact these smooth gradients have on optimization routines in the presence of discontinuous subgradients. 
\subsection{Collision Detection}
    \begin{figure}[t!]
        \centering
        \includegraphics[width=.8\linewidth]{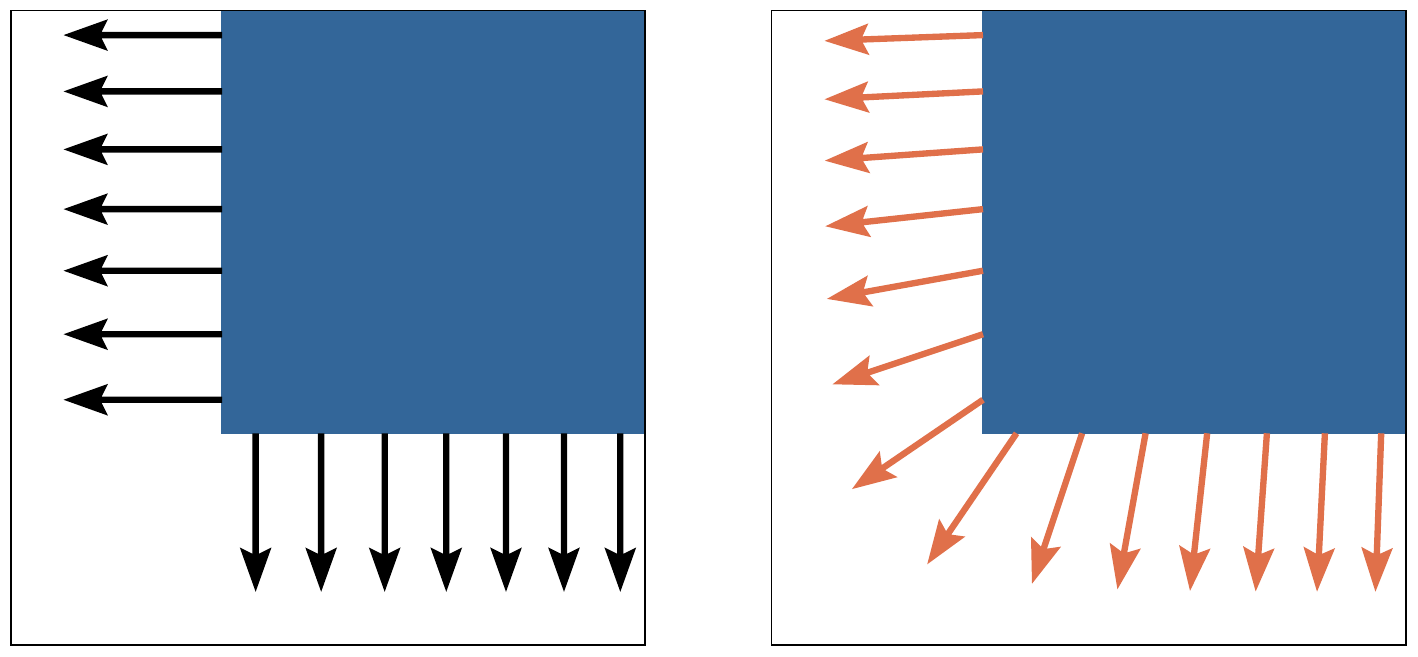}

        \definecolor{CUST}{HTML}{E0704A}
        \begin{tikzpicture}

            \draw[black, line width=2pt] (-3.7, 1.77) -- (-3.0, 1.77);
            \node[anchor=west] at (-2.8, 1.8) {$\kappa=0$};

            \draw[CUST, line width=2pt] (1.3, 1.77) -- (2.0, 1.77);
            \node[anchor=west] at (2.2, 1.8) {$\kappa=0.01$};
        \end{tikzpicture}
        \caption{Contact normal vectors from an optimization-based differentiable collision detection routine with and without relaxed differentiation. With no relaxation ($\kappa=0$), the direction of the contact normal switches immediately as the closest point moves from one face to another. When the relaxed gradients are used ($\kappa=0.01$), the contact normal smoothly transitions between the faces.}
        \label{fig:dcol}
    \end{figure}
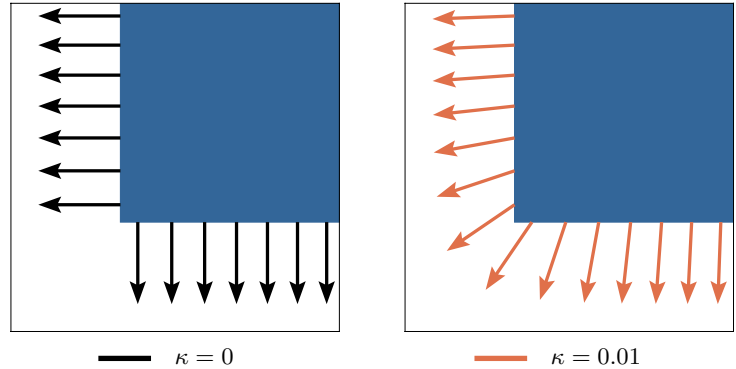
Collision detection between convex shapes can be formulated as a convex optimization problem, both in terms of the closest point between shapes \cite{gilbert1988}, and in terms of the minimum scale factor \cite{tracy2023b}. For the former, we introduce two points in a world frame $p_i \in \R{3}$, and two polytopes described with $A_i p_i \leq b_i$. By constraining each point to be within a polytope, a QP is used to solve for the closest point between these two shapes,
\begin{mini}
    {p_1, p_2}{ \|p_1 - p_2\|}{\label{gjk}}{}
    \addConstraint{A_1 p_1}{\leq b_1}
    \addConstraint{A_2 p_2}{\leq b_2.}
\end{mini}
Using a differentiable QP solver, the gradient of the closest distance between shapes with respect to the positions of the polytopes can be computed resulting in the contact normal vectors. 

In this example, we examine collision detection between two squares and the behavior of these contact normals in the presence of sharp corners. As shown in Fig. \ref{fig:dcol}, the contact normals are evaluated at a strict $\kappa=0$ and a relaxed $\kappa=0.01$. In the case of $\kappa=0$, the contact normals are (correctly) exactly normal to the surface, and as soon as the closest point shifts from one face to the other, the contact normals immediately rotate $90^\circ$. While this is expected behavior, this discontinuity in the gradient can prove troublesome for simulation and control algorithms that rely on these contact normals not changing too quickly. Alternatively, with a relaxed $\kappa=0.01$, the contact normal smoothly rotates the 90 degrees as the face of the closest point changes. This is a result of the logarithmic barrier smoothing out the sharp corner, and allows for continuous and smooth gradients even in the presence of the discontinuity. 
\section{Conclusions}
In this paper, we outline shortcomings with existing differentiable optimization tools, namely the nonsmoothness of the gradients near inequality constraints and the inability to handle infeasible problems, and propose solutions to both of these problems. By relaxing the solution to an optimization problem from tight tolerances to an intentionally relaxed logarithmic barrier, unique and smooth gradients can be computed even from sharp edges in the feasible set. This relaxation is straightforward and leverages existing routines within existing primal-dual interior-point solvers. We also introduce an always-feasible quadratic program where hard constraints are converted into $\ell_1$-penalties, and devise a customized algorithm for solving problems of this form with limited added computational overhead. With both of these innovations, consistent and reliable smooth gradients are demonstrated in common robotic tasks where smoothness is a priority. Our fully differentiable and parallelizable solver written in JAX is available at \url{www.github.com/kevin-tracy/qpax}.
\begin{credits}
\subsubsection{\discintname}
The authors have no competing interests to declare that are
relevant to the content of this article..
\end{credits}
%
%
%
\bibliographystyle{splncs04}
\bibliography{references}
%




\end{document}